# POISSON CALCULUS FOR SPATIAL NEUTRAL TO THE RIGHT PROCESSES[1]

By Lancelot F. James

*Hong Kong University of Science and Technology*

Neutral to the right (NTR) processes were introduced by Doksum in 1974 as Bayesian priors on the class of distributions on the real line. Since that time there have been numerous applications to models that arise in survival analysis subject to possible right censoring. However, unlike the Dirichlet process, the larger class of NTR processes has not been used in a wider range of more complex statistical applications. Here, to circumvent some of these limitations, we describe a natural extension of NTR processes to arbitrary Polish spaces, which we call spatial neutral to the right processes. Our construction also leads to a new rich class of random probability measures, which we call NTR species sampling models. We show that this class contains the important two parameter extension of the Dirichlet process. We provide a posterior analysis, which yields tractable NTR analogues of the Blackwell–MacQueen distribution. Our analysis turns out to be closely related to the study of regenerative composition structures. A new computational scheme, which is an ordered variant of the general Chinese restaurant processes, is developed. This can be used to approximate complex posterior quantities. We also discuss some relationships to results that appear outside of Bayesian nonparametrics.

**1. Introduction.** Doksum [9] considered a nonparametric Bayesian analysis based on neutral to the right (NTR) priors. These priors are random probability measures defined on the real line, $\mathcal{R}$, that include the popular Dirichlet process (see [13] and [16]). Within Bayesian nonparametric statistics, the NTR process serves as one of the important classes of models. In particular, there have been numerous applications to models that

Received May 2003; revised March 2005.
[1]Supported in part by RGC Grant HKUST-6159/02P and DAG 01/02.BM43 of HK-SAR.
*AMS 2000 subject classifications.* Primary 62G05; secondary 62F15.
*Key words and phrases.* Bayesian nonparametrics, inhomogeneous Poisson process, Lévy processes, neutral to the right processes, regenerative compositions, survival analysis.

arise in survival analysis subject to possible right censoring. On the other hand, unlike the Dirichlet process, the larger class of NTR processes has not been used in a wider range of statistical applications. That is, for instance, there are no general NTR analogues of the important class of (kernel based) Dirichlet process mixture models. See, for example, [32] and [23] for further background and references on Dirichlet process mixture models.

A goal of this article is to begin to answer the question of how one can possibly use NTR processes in a wider context, as has been the case for the Dirichlet process. One of the limitations of NTR processes is that they are only defined on the real line. The other limitation, which is perhaps more severe, is that as of yet we do not have *tractable* NTR analogues of the Blackwell–MacQueen [3] Pólya urn distribution associated with the Dirichlet process. The Blackwell–MacQueen distribution is well known to be the exchangeable distribution derived from a Dirichlet process, and its theoretical understanding and practical implementation are crucial in complex models. To circumvent some of these limitations, we describe a natural extension of NTR processes defined on an arbitrary Polish space $\mathscr{S} = \mathcal{R}^+ \times \mathscr{X}$, which we call spatial NTR processes. Here $\mathcal{R}^+$ denotes the positive real line and $\mathscr{X}$ is an arbitrary Polish space. Our construction also leads to a rich class of random probability measures on $\mathscr{X}$, which we call *NTR species sampling*. We provide a detailed analysis of these models and obtain properties analogous to the Dirichlet process. In particular, we provide a description of the posterior distribution of spatial NTR processes and, more importantly, we give a detailed analysis of the NTR analogues of the Blackwell–MacQueen distribution.

Such an analysis parallels, in part, the results of Antoniak [1] (see also [12]) and Lo [32] for the Dirichlet process. These works involve characterizations based on random partitions of the integers $\{1, \ldots, n\}$ and were derived using nontrivial combinatorial arguments. The structure of general NTR processes is more complex than that of the Dirichlet process and an approach using direct combinatorial analysis is considerably more challenging. We circumvent such issues by applying the *Poisson process partition calculus* discussed by James [24, 26]. This also paves the way for a straightforward derivation of the posterior distribution of spatial NTR processes. Using these results, we develop a new computational scheme related to the general Chinese restaurant process (see [37], page 60 and [23]), which now allows one to sample from the exchangeable distributions derived from NTR processes.

It is important to note that although Bayesian applications of NTR processes to complex statistical models have been limited, the use of these processes appears often in other important contexts. Doksum ([9], Theorem 3.1) showed that one can describe an NTR distribution function $F$ on $\mathcal{R}^+$ via positive Lévy processes, $Z$, on $\mathcal{R}^+$ as

$$1 - F(t) = S(t) = e^{-Z(t)}, \tag{1}$$

where $S$ denotes the survival distribution of a random variable $T$ from $F$. The Lévy process $Z$ is an increasing independent increment process that satisfies $Z(0) = 0$ and $\lim_{t \to \infty} Z(t) = \infty$ a.s. That is, $T|F$ has survival distribution $P(T > t|F) = e^{-Z(t)}$. Importantly, the representation in (1) shows that NTR survival processes essentially coincide with the class of exponential functionals of possibly inhomogeneous, nonnegative Lévy processes. Such objects and more general exponential functionals of Lévy processes, such as Brownian motion, have been extensively studied by probabilists with applications, for instance, to finance. The NTR models also arise in coalescent theory, which has applications in genetics and physics, as seen, for example, in [36], Proposition 26. See also [2] and [5]. Noting some of these connections, Epifani, Lijoi and Prünster [11] applied techniques from those manuscripts to obtain expressions for the moments of mean functionals of NTR models and, as we also do here, highlighted some of the connections to these areas outside of Bayesian nonparametric statistics. The mean functional can be described explicitly as

$$(2) \qquad I = \int_0^\infty tF(dt) = \int_0^\infty S(t)\,dt = \int_0^\infty e^{-Z(t)}\,dt.$$

It is a significant object, which has interesting interpretations in a variety of fields. We describe how this process is related to the study of the Blackwell–MacQueen analogue derived from NTR processes. Moreover, we discuss how our work is closely related to the recent work of Gnedin and Pitman [18] on *regenerative composition structures*.

**2. Spatial neutral to the right processes.** Suppose that $(T, X)$ are random elements on the Polish space $\mathscr{S}$ that have distribution $F(ds, dx)$ for $(s, x) \in \mathscr{S}$. Here we would like to extend the definition of an NTR process to model $F(ds, dx)$ as a random probability measure such that its marginal $F(ds, \mathscr{X})$ is an NTR process. While the representation in (1) is quite useful for calculations, it is not immediately obvious how one can use this definition to extend an NTR process to $\mathscr{S}$. The known exception is the Dirichlet process that can be defined on arbitrary spaces. To do this, we first recall that if $F$ is an NTR process on $\mathcal{R}^+$, then its cumulative hazard $\Lambda(ds) := F(ds)/S(s-)$ is a nonnegative Lévy process; in other words, $\Lambda$ is a completely random measure (see [30]). This observation and alternative idea for modeling via cumulative hazards is due to the important work of Hjort [21]. Note that $Z$ in (1) is also a completely random measure. Moreover, an important aspect of our results relies on the fact that there is one-to-one distributional correspondence between a particular $Z$ and $\Lambda$. Specifically, if $J_j$ represents a random jump of $\Lambda$ taking its values in $[0, 1]$, then $-\log(1 - J_j)$ is the jump of a corresponding $Z$ taking its values in $\mathcal{R}^+$. Hence if we initially model $Z$ and $\Lambda$ as completely random measures without a drift component and fixed points of discontinuity, they may both be represented as linear functionals of a common Poisson random measure. Importantly, one then may give precise meaning to the distributional equivalences $\mathbb{P}(T \in ds|F) = F(ds) := S(s-)\Lambda(ds) := e^{-Z(s-)}\Lambda(ds)$, where $F$ is an NTR process.

Our construction now proceeds by extending $\Lambda$ and $Z$ to completely random measures on $\mathscr{S}$, using a representation in terms of a Poisson random measure. Let $N$ denote a Poisson random measure on some Polish space $\mathscr{W} = [0,1] \times \mathscr{S}$ with mean intensity

$$\mathbb{E}[N(du, ds, dx)|\nu] = \nu(du, ds, dx) := \rho(du|s)\Lambda_0(ds, dx).$$

Here $\rho$ is a Lévy density that will determine the conditional distribution of the jumps of $\Lambda$ and $Z$. Furthermore, without loss of generality, we assume that $\int_0^1 u\rho(du) = 1$ and, hence, $u\rho(du)$ is a probability density function. The intensity $\nu$ is chosen such that $\Lambda_0(ds, dx) := F_0(ds, dx)/S_0(s-)$ is by definition a hazard measure on $\mathscr{S}$, where $F_0$ represents a prior specification for the distribution of $F$ on $\mathscr{S}$, and $S_0$ is the corresponding survival function on $\mathcal{R}^+$. See [31], A5.3, for formal details of hazard measures on abstract spaces. Note that in [28], Proposition 25.28, the hazard measure is also called a *natural compensator* of a random measure defined as $\delta_{T,X}$. We denote the Poisson law of $N$ with intensity $\nu$ as $\mathbb{P}(dN|\nu)$. The Laplace functional for $N$, which plays an important role in our analysis, is defined as

$$\mathbb{E}[e^{-N(f)}|\nu] = \int_{\mathbb{M}} e^{-N(f)} \mathbb{P}(dN|\nu) = e^{-\mathscr{G}(f)},$$

where for any positive $f$, $N(f) = \int_{\mathscr{W}} f(x)N(dx)$ and $\mathscr{G}(f) = \int_{\mathscr{W}} (1 - e^{-f(x)}) \times \nu(dx)$, and $\mathbb{M}$ denotes the space of boundedly finite measures on $\mathscr{W}$ (see [6]). A measure, say $N$, is boundedly finite if for each bounded set $A$, $N(A) < \infty$. See also [28], Chapter 12, for a discussion of Poisson random measures and the unicity property of Laplace functionals.

Now the specifications above imply that $\Lambda(ds, dx) := \int_0^1 uN(du, ds, dx)$ is a *completely random hazard measure* on $\mathscr{S}$ with mean $E[\Lambda(ds, dx)] = \Lambda_0(ds, dx)$ and there is a corresponding $Z(ds, dx) := \int_0^1 [-\log(1-u)] \times N(du, ds, dx)$. In particular, $-\log S(t-) := Z(t-) = \int_{\mathscr{W}} [-I\{s < t\}\log(1-u)]N(du, ds, dx)$. Now using these facts we define a *spatial neutral to the right* (SPNTR) random probability measure on $\mathscr{S}$ as

(3) $\qquad \mathbb{P}(T \in dt, X \in dx|F) := F(dt, dx) = S(t-)\Lambda(dt, dx).$

Defining $\Lambda(ds) := \Lambda(ds, \mathscr{X})$, it follows that $F(ds) := S(s-)\Lambda(ds)$ is an NTR process and, furthermore, $\mathbb{E}[F(dt, dx)] = S_0(t-)\Lambda_0(dt, dx) = F_0(dt, dx)$. See Section 5 for more details.

REMARK 1. The choice of

$$(4) \qquad \rho(du|s)\Lambda_0(ds,dx) = c(s)u^{-1}(1-u)^{c(s)-1}\,du\,\Lambda_0(ds,dx)$$

for $c(s)$ a positive function yields a natural extension of Hjort's [21] beta cumulative hazard process to beta processes on $\mathscr{S}$. Equivalently, this specification defines beta-Stacy or beta-neutral distribution functions on $\mathscr{S}$. See [33, 41] and [21], Section 7A, for such processes defined on $\mathcal{R}^+$. The case of the Dirichlet process with shape parameter $\theta F_0$ is obtained by choosing $c(s) = \theta S_0(s-)$. Our construction of spatial NTR processes is influenced by the work of James and Kwon [27], who first gave an explicit construction of spatial beta-neutral processes on $\mathscr{S}$ via ratios of two independent gamma processes.

REMARK 2. Given the specifications in (3), we extend this definition to include prior fixed points of discontinuity $\{(s_1, w_1), \ldots, (s_k, w_k)\}$ in $\mathscr{S}$ as

$$(5) \qquad \tilde{F}_k(ds,dx) = \left[e^{-Z(s-)} \prod_{\{l\,:\,s_l<s\}} (1-U_l)\right]\tilde{\Lambda}_k(ds,dx),$$

where $\tilde{\Lambda}_k(ds,dx) = \Lambda(ds,dx) + \sum_{l=1}^k U_l \delta_{s_l,w_l}(ds,dx)$ is defined such that independent of $\Lambda$, $U_j$ are independent random variables on $[0,1]$ with distribution $H_j$ for $j = 1, \ldots, k$. We call $\tilde{F}_k$ a general spatial NTR process.

REMARK 3. The log mapping that we use can be deduced, for instance, from [7] and [8], Proposition 2. This type of correspondence is actually noted, albeit less explicitly, in [21] and is also used in related contexts without specific mention of NTR processes; see, for instance, [36], Proposition 26. In particular, if $\tau$ is a Lévy measure that specifies the conditional distribution of the jumps of $Z$, then by writing $\tau(dy|s) := \tau(y|s)\,dy$ and $\rho(du|s) := \rho(u|s)\,du$, the relationship between the Lévy measures of $Z$ and $\Lambda$ is described by

$$\tau(y|s) = e^{-y}\rho(1-e^{-y}|s) \qquad \text{for } y \in \mathcal{R}^+ \text{ or}$$
$$\rho(u|s) = (1-u)^{-1}\tau(-\log(1-u)|s) \qquad \text{for } u \in [0,1].$$

Note that if $\rho(du|s) := \rho(du)$, then we say that the relevant processes are homogeneous.

**3. Posterior analysis.** Similar to the case of the Dirichlet process, we consider the following setup. Suppose that $(T_i, X_i)|F$ are i.i.d. pairs with common distribution $F$ for $i = 1, \ldots, n$ and suppose the law of $F$ is modeled as a spatial NTR process. This description yields a joint distribution of $(\mathbf{T}, \mathbf{X}) = \{(T_1, X_1), \ldots, (T_n, X_n)\}$ and $F$. We are interested in the Bayesian disintegration of this joint distribution in terms of the posterior distribution of $F|\mathbf{T}, \mathbf{X}$ and the marginal distribution of $(\mathbf{T}, \mathbf{X})$. Since $\Lambda$, $Z$ and $F$ are all functionals of $N$, we work instead with the joint distribution of $(\mathbf{T}, \mathbf{X}, N)$,

$$
\begin{aligned}
&\left[\prod_{i=1}^{n} S(T_i-)\Lambda(dT_i, dX_i)\right]\mathbb{P}(dN|\nu) \\
&\quad = \pi(dN|\mathbf{T}, \mathbf{X})\mathscr{M}(dT_1, dX_1, \ldots, dT_n, dX_n),
\end{aligned}
\tag{6}
$$

where $\pi(dN|\mathbf{T}, \mathbf{X})$ denotes the desired posterior distribution of $N|\mathbf{T}, \mathbf{X}$ and

$$
\begin{aligned}
\mathscr{M}(d\mathbf{T}, d\mathbf{X}) &= \mathscr{M}(dT_1, dX_1, \ldots, dT_n, dX_n) \\
&= \int_{\mathcal{M}} \left[\prod_{i=1}^{n} F(dT_i, dX_i)\right] \mathcal{P}(dN|\nu)
\end{aligned}
\tag{7}
$$

is the important exchangeable marginal distribution of $(\mathbf{T}, \mathbf{X})$. The $\mathscr{M}$ denotes the general analogue of the Blackwell–MacQueen Pólya urn, and hence is crucial to both theoretical understanding and practical implementation of the general class of spatial NTR processes. We will describe the posterior distribution given $(\mathbf{T}, \mathbf{X})$ in Section 4, and give a detailed analysis of $\mathscr{M}$ and related quantities in Section 5. We first explain some key elements of the analysis.

3.1. *The role of random partitions and order statistics.* It is clear that one can always represent $(\mathbf{T}, \mathbf{X})$ as $(\mathbf{T}^*, \mathbf{X}^*, \mathbf{p})$, where, using notation similar to Lo [32], $(\mathbf{T}^*, \mathbf{X}^*) = \{(T_1^*, X_1^*), \ldots, (T_{n(\mathbf{p})}^*, X_{n(\mathbf{p})}^*)\}$ denotes the distinct pairs of observations within the sample and where $\mathbf{p} = \{E_1, \ldots, E_{n(\mathbf{p})}\}$ stands for a partition of $\{1, \ldots, n\}$ of size $n(\mathbf{p}) \le n$ that records which observations within the sample are equal. The number of elements in the $j$th cell, $E_j := \{i : (T_i, X_i) = (T_j^*, X_j^*)\}$, of the partition is indicated by $e_j$ for $j = 1, \ldots, n(\mathbf{p})$, so that $\sum_{j=1}^{n(\mathbf{p})} e_j = n$. It follows that the marginal distribution of $(\mathbf{T}, \mathbf{X})$, say $\mathscr{M}$, can be expressed in terms of a conditional distribution of $\mathbf{T}, \mathbf{X}|\mathbf{p}$, which is the same as a conditional distribution of the unique values $\mathbf{T}^*, \mathbf{X}^*|\mathbf{p}$, and the marginal distribution of $\mathbf{p}$. The marginal distribution of $\mathbf{p}$, denoted as $\pi(\mathbf{p})$ or $p(e_1, \ldots, e_{n(\mathbf{p})})$, is an *exchangeable partition probability function* (EPPF), that is, a probability distribution on $\mathbf{p}$ which is exchangeable in its arguments and depends only on the size of each cell. The best known case of an EPPF is the variant of the Ewens sampling formula (ESF) (see [1, 12]) associated with the Dirichlet process with total mass $\theta$, given as

$$\frac{\theta^{n(\mathbf{p})}\Gamma(\theta)}{\Gamma(\theta + n)} \prod_{j=1}^{n(\mathbf{p})} \Gamma(e_j).$$

Additionally, since a Dirichlet process is a special case of what are called species sampling models, the distribution of $\mathbf{T}, \mathbf{X}|\mathbf{p}$ is such that the unique pairs $(T_j^*, X_j^*)$ are i.i.d. with distributions $F_0$. We note that the marginal distribution and, naturally, the posterior distribution of the Dirichlet process depend only on the counts $e_j$ and the unique values. The structure of $\mathscr{M}$ for general NTR processes is considerably more complex. However, as we will explain, what is interesting is that they do have a natural interpretation in terms of classical survival models. One can think of $\mathbf{T}^*$ as the collection of the *unordered* distinct times to death of individuals in a sample of size $n$. In this sense, the count $e_j$ represents the number of deaths at time $T_j^*$. Additionally, it is well known that the posterior distribution of NTR processes also depends on the *number at risk* at a given time, say $t$, which can be defined as $Y_n(t) = \sum_{i=1}^n I\{T_i > t\}$. We have discovered that to simplify the expressions for $\mathscr{M}$, it is necessary not only to know the number of deaths, but also to know the number at risk at the unique times. Instead of working directly with $\mathbf{T}^*$, we do this by using its ordered values.

That is, let $T_{(1:n)} > T_{(2:n)} > \cdots > T_{(n(\mathbf{p}):n)}$ denote an ordering of the unique values $\{T_1^*, \ldots, T_{n(\mathbf{p})}^*\}$. This collection represents the *ordered* unique times of death. Note that we work with the pairs $(T_{(j:n)}, X_j^*)$, where $X_j^*$ is simply the unique value treated as the concomitant of $T_{(j:n)}$. That is, we do not order the $\mathbf{X}$ values; in fact, some spaces $\mathscr{X}$ do not have a natural ordering. Associated with this, let $\mathbf{m} = \{E_{(1)}, \ldots, E_{(n(\mathbf{p}))}\}$ denote the collection of sets $E_{(j)} = \{i : T_i = T_{(j:n)}\}$ for $j = 1, \ldots, n(\mathbf{p})$. That is, $E_{(j)}$ is the collection of values equal to the $j$th largest unique death time. Similar to $e_j$, let $m_j = |E_{(j)}|$ denote the number of deaths at the $j$th largest unique death time, $T_{(j:n)}$, for $j = 1, \ldots, n(\mathbf{p})$. There are of course $n(\mathbf{p})!$ possible orderings of $\mathbf{T}^*$. This implies that given a partition $\mathbf{p} = \{E_1, \ldots, E_{n(\mathbf{p})}\}$, the collection $\{m_1, \ldots, m_{n(\mathbf{p})}\}$ [resp. $(\mathbf{m})$] takes its values over the symmetric group, say $S_{n(\mathbf{p})}$, of all $n(\mathbf{p})!$ permutations of $\{e_1, \ldots, e_{n(\mathbf{p})}\}$ [of $(\{E_1, \ldots, E_{n(\mathbf{p})}\})$]. Notice now that, for each $s$,

$$Y_n(s) = \sum_{i=1}^n I\{T_i > s\} = \sum_{j=1}^{n(\mathbf{p})} e_j I\{T_j^* > s\} = \sum_{l=1}^{n(\mathbf{p})} m_l I\{T_{(l:n)} > s\}.$$

Hence for $j = 1, \ldots, n(\mathbf{p})$, we can define $r_{j-1} := Y_n(T_{(j:n)}) = \sum_{l=1}^{j-1} m_l$, which denotes the number larger than the $j$th largest unique value. Note that $r_0 = 0$ and $r_{n(\mathbf{p})} = n$; additionally, $r_j = r_{j-1} + m_j$. What is important is that, unlike $\mathbf{p}$, the collection $\{E_{(1)}, \ldots, E_{(n(\mathbf{p}))}\}$ completely determines $(r_j)$ via the $(m_j)$; that is, $\mathbf{m}$ contains the relevant information in $\mathbf{p}$. We will often refer to $(\mathbf{m}, \mathbf{p})$ rather than $\mathbf{m}$ to remind the reader of the dependence of $\mathbf{m}$ on $\mathbf{p}$.

REMARK 4. See [34, 37, 38] for a general overview of the EPPF concept and see [23, 24, 26] for its relevance to general marginal exchangeable distributions that arise in a Bayesian context.

REMARK 5. One of the earliest applications of the Ewens sampling formula is in population genetics. It is quite interesting to note that, as described by Donnelly and Joyce ([10], page 230), one may also interpret the $(T_{(j:n)})$ as the ordering of genetic types (alleles) of individuals, where new alleles arise by mutation and the alleles present in the population or in a sample at a given time may be ordered by age. An interesting by-product of our work is that it actually yields the explicit distribution for large classes of such models. A simple description will be given in Proposition 5.2.

**4. The posterior distribution of spatial NTR processes.** In this section we describe formally the posterior distribution of spatial NTR processes given the data $(\mathbf{T}, \mathbf{X})$. Note here that we will characterize the posterior via the ordered values rather than $\mathbf{T}^*$. Since we are conditioning on $(\mathbf{T}, \mathbf{X})$, these are equivalent notions. We first describe the result for no fixed points of discontinuity and then discuss how one easily obtains the extension in Section 4.1. The proof is delayed until the Appendix.

PROPOSITION 4.1. *Let F be a spatial NTR process defined by the Poisson random measure N with mean intensity $\nu(du, ds, dx) = \rho(du|s)\Lambda_0(ds, dx)$; $\Lambda$ is its corresponding Lévy hazard measure. Suppose that $(T_i, X_i)|F$ are i.i.d. F for $i = 1, \ldots, n$. Then:*

(i) *The posterior distribution of $N|\mathbf{T}, \mathbf{X}$ is equivalent to the distribution of the random measure $N_n^* = N_n + \sum_{j=1}^{n(\mathbf{p})} \delta_{J_{j,n}, T_{(j:n)}, X_j^*}$, where, conditional on $(\mathbf{T}, \mathbf{X})$, $N_n$ is a Poisson random measure with intensity*

(8) $$\nu_n(du, ds, dx) = (1-u)^{Y_n(s)} \rho(du|s) \Lambda_0(ds, dx).$$

*Additionally, the $(J_{j,n})$ are conditionally independent of $N_n$ and are mutually independent with distributions specified by*

$$\mathbb{P}(J_{j,n} \in du | T_{(j:n)}) \propto u^{m_j}(1-u)^{r_j - 1} \rho(du | T_{(j:n)})$$

*for $j = 1, \ldots, n(\mathbf{p})$.*

(ii) *The posterior distribution of $\Lambda$ given $(\mathbf{T}, \mathbf{X})$ is equivalent to the law of the Lévy hazard measure,*

$$\Lambda_n^*(ds, dx) = \int_0^1 u N_n^*(du, ds, dx)$$

$$= \Lambda_n(ds, dx) + \sum_{j=1}^{n(\mathbf{p})} J_{j,n} \delta_{T_{(j:n)}, X_j^*}(ds, dx),$$

where $\Lambda_n(ds, dx) = \int_0^1 u N_n(du, ds, dx)$ is a Lévy hazard measure with Lévy measure as in (8) and where the $(J_{j,n})$ are conditionally independent of $\Lambda_n$.

(iii) *The posterior distribution of the corresponding $Z$ process is equivalent to the the law of the random measure*

$$Z_n^*(ds, dx) = Z_n(ds, dx) + \sum_{j=1}^{n(\mathbf{p})} Z_{j,n} \delta_{T_{(j:n)}, X_j^*}(ds, dx),$$

where $Z_n(ds, dx) = \int_0^1 [-\log(1-u)] N_n(du, ds, dx)$ and each $Z_{j,n} = -\log(1 - J_{j,n})$ with distribution

$$\mathbb{P}(Z_{j,n} \in dy | T_{(j:n)}) := H_j^*(d(1 - e^{-y})) \propto (1 - e^{-y})^{m_j} e^{-r_j - 1} \tau(dy | T_{(j:n)}).$$

(iv) *Additionally, the posterior distribution of $F$ is equivalent to the conditional law, given $(\mathbf{T}, \mathbf{X})$, of the random probability measure $F_n^*(ds, dx)$ expressed as*

$$e^{-Z_n(s-)} \left[ \prod_{\{j : T_{(j:n)} < s\}}^{n(\mathbf{p})} (1 - J_{j,n}) \right] \Lambda_n(ds, dx) + \sum_{j=1}^{n(\mathbf{p})} \tilde{P}_{j:n} \delta_{T_{(j:n)}, X_j^*}(ds, dx),$$

where $\tilde{P}_{j,n} = e^{-Z_n(T_{(j:n)}-)} J_{j,n} \prod_{l=j+1}^{n(\mathbf{p})} (1 - J_{l,n})$. It follows that the Bayesian prediction rule is given by $E[F_n^*(ds, dx) | \mathbf{T}, \mathbf{X}]$, which can be expressed in several ways.

REMARK 6. Note that due to symmetry, one has the equivalence in distribution of

$$\sum_{j=1}^{n(\mathbf{p})} J_{j,n} \delta_{T_{(j:n)}, X_j^*}(ds, dx) = \sum_{j=1}^{n(\mathbf{p})} J_{j,n}^* \delta_{T_j^*, X_j^*}(ds, dx),$$

where the random variables $(J_{j,n}^*)$ are mutually independent with marginal distributions $\mathbb{P}(J_{j,n}^* \in ds | T_j^*) \propto u^{e_j}(1-u)^{Y_n(T_j^*)} \rho(du | T_j^*)$. Recall that $Y_n(T_{(j:n)}) = r_{j-1}$.

4.1. *Remarks on prior fixed points of discontinuity.* We have so far omitted any discussion on the form of the posterior distribution when there are prior points of discontinuity as in $\tilde{\Lambda}_k$ defined in (5). In fact, the analysis is essentially already contained in our results. Recall that for $n \geq 1$ the posterior process for $\Lambda$ in the complete data is $\Lambda_n^* = \Lambda_n + \sum_{j=1}^{n(\mathbf{p})} J_{j,n} \delta_{T_{(j:n)}, X_j^*}$, where the $(J_{j,n})$ are conditionally independent of $\Lambda_n$. Using the fact that $\tilde{\Lambda}_k$ and $\Lambda_n^*$ are the same structurally, one can simply let $n(\mathbf{p})$ play the role of $k$ and let $\{U_l, s_l, w_l\}$ play the role of $\{J_{j,n}, T_{(j:n)}, X_j^*\}$. Let $n_l = |\{i : (T_i, X_i) = (s_l, w_l)\}|$ for $l = 1, \ldots, k$. In addition, let $(T_{(j:n)}, X_j^*)$ denote $0 \leq n(\mathbf{p}) \leq n$ unique values distinct from $\{(s_1, w_1), \ldots, (s_k, w_k)\}$. Then it is easy to see that the posterior distribution of $\tilde{\Lambda}_k$ is of the form

$$\tilde{\Lambda}_{n,k}^* = \Lambda_n + \sum_{l=1}^{k} U_{l:n} \delta_{s_l, w_l} + \sum_{j=1}^{n(\mathbf{p})} J_{j,n} \delta_{T_{(j:n)}, X_j^*},$$

where $\mathbb{P}\{U_{l:n} \in du|s_l\} \propto u^{n_l}(1-u)^{Y_n(s_l)} H_l(du)$ for $l = 1, \ldots, k$. Note here we use $Y_n(s) = \sum_{i=1}^{n} I\{T_i > s\}$.

REMARK 7. Note that marginalizing over $\mathscr{X}$, the result in Proposition 4.1 reduces to the appropriate analogous results for NTR processes described in [9, 14, 15, 21, 29]. However, we shall present a considerably streamlined and direct proof that uses a methodology applicable to a much wider class of random probability measures on abstract spaces. Note moreover that there is no analogue of Proposition 4.1(i) appearing in those works. The distribution of $F_n^*(\infty, dx)$ corresponds to the posterior distribution of a new class of random probability measures, which we discuss in more detail in Section 5.3.

**5. Analysis of NTR generalizations of the Blackwell–MacQueen distribution.** We now present a detailed analysis of the marginal distribution $\mathscr{M}$ and related quantities. We give details for Lemmas 5.1 and 5.2 in the Appendix. First we introduce some additional notation. For a homogeneous $\rho$ or $\tau$ and for $\omega \geq 0$, let

$$\phi(\omega) = \int_0^\infty (1 - e^{-\omega y}) \tau(dy) = \int_0^1 (1 - (1-u)^\omega) \rho(du)$$
$$= \int_0^1 \omega (1-u)^{\omega-1} \left[ \int_u^1 \rho(dv) \right] du.$$

This is the *Lévy exponent* defined by the Laplace transform of a homogeneous $Z$ process. For integers $(i, k)$, let

$$\psi_{i,k}(s) = \int_0^\infty (1 - e^{-yi}) e^{-yk} \tau(dy|s) = \int_0^1 (1 - (1-u)^i)(1-u)^k \rho(du|s).$$

In the homogeneous case, set $\psi_{i,k} = \int_0^\infty (1 - e^{-yi}) e^{-yk} \tau(dy)$ and note that for each $j$, $\phi(j) = \psi_{j,0} = \int_0^\infty (1 - e^{-jy}) \tau(dy)$. Finally, we define cumulants

$$\kappa_{m_j, r_j-1}(\rho|s) = \int_0^1 u^{m_j}(1-u)^{r_j-1} \rho(du|s)$$

and

$$\kappa_{m_j, r_j-1}(\rho) = \int_0^1 u^{m_j}(1-u)^{r_j-1} \rho(du).$$

Our first task will be to obtain a nice expression for the expectation of the product of survival functions that appears in (6). First notice that

$$(9) \quad \left[\prod_{i=1}^{n} S(T_i-)\right] = \left[\prod_{j=1}^{n(\mathbf{p})} S(T_j^*-)^{e_{j,n}}\right] = \left[\prod_{j=1}^{n(\mathbf{p})} S(T_{(j:n)}-)^{m_j}\right].$$

These equivalences lead to the following result.

LEMMA 5.1. *Let $\nu(du, ds, dx) = \rho(du|s)\Lambda_0(ds, dx)$ be the mean intensity of a Poisson random measure $N$. Then*

$$\mathbb{E}\left[\prod_{j=1}^{n(\mathbf{p})} S(T_{(j:n)}-)^{m_j} | \nu\right] = \prod_{j=1}^{n(\mathbf{p})} e^{-\int_0^{T_{(j:n)}} \psi_{m_j, r_{j-1}}(s)\Lambda_0(ds)}.$$

*The expression reduces to $\prod_{j=1}^{n} \exp(-\int_0^{T_{(j:n)}} \psi_{1,j-1}(s)\Lambda_0(ds))$ when there are no ties.*

Lemma 5.1 is instrumental in obtaining the following initial description of $\mathscr{M}$.

LEMMA 5.2. *Let $\mathscr{M}$ denote the exchangeable distribution of $(\mathbf{T}, \mathbf{X})$ defined in (7). Then $\mathscr{M}(d\mathbf{T}, d\mathbf{X})$ can be expressed as*

$$\left[\prod_{j=1}^{n(\mathbf{p})} e^{-\int_0^{T_{(j:n)}} \psi_{m_j, r_{j-1}}(s)\Lambda_0(ds)} \kappa_{m_j, r_{j-1}}(\rho|T_{(j:n)})\right] \prod_{l=1}^{n(\mathbf{p})} \Lambda_0(dT_l^*, dX_l^*).$$

We now show how one can obtain calculations using $\mathscr{M}$. For each $\mathbf{m} \in S_{n(\mathbf{p})}$ and integrable function $g(\mathbf{T})$, define

$$L(g; \mathbf{m}) = \int_0^\infty \int_{t_{n(\mathbf{p})}}^\infty \cdots \int_{t_2}^\infty g((\mathbf{t}, \mathbf{m})) \prod_{j=1}^{n(\mathbf{p})} e^{-\int_0^{t_j} \psi_{m_j, r_{j-1}}(s)\Lambda_0(ds)}$$
$$\times \kappa_{m_j, r_{j-1}}(\rho|t_j)\Lambda_0(dt_j),$$

where $t_1 > t_2 > \cdots > t_{n(\mathbf{p})}$ denotes one of $n(\mathbf{p})!$ orderings of the unique values. With some abuse of notation, the vector $(\mathbf{t}, \mathbf{m}) = (\mathbf{t})$ denotes the collection of $n$ points whose $n(\mathbf{p})$ unique values are ordered according to $\mathbf{m}$. For example, suppose one has the function $g(T_1, T_2, T_3)$. Then in the instance where $T_1 = T_2 < T_3$, one has $n(\mathbf{p}) = 2$ unique values and one evaluates $g(T_{(2:2)}, T_{(2:2)}, T_{(1:2)})$ or, using the notation above, $g(t_2, t_2, t_1)$.

We now use $L$ to obtain very general formulae for expected values of complex integrals of NTR processes. This plays a key role in obtaining the EPPF $\pi(\mathbf{p})$ and related quantities.

LEMMA 5.3.  *Assume that the random functional $\mathscr{I}(g) = \int g(\mathbf{t}) \prod_{i=1}^n F(dt_i)$ is integrable, where $F$ is an NTR process specified by the Poisson law $\mathbb{P}(dN|\nu)$. Then it follows from Lemma 5.2 that*

$$\mathbb{E}[\mathscr{I}(g)|\nu] = \sum_{\mathbf{p}} \left[ \sum_{\mathbf{m} \in S_{n(\mathbf{p})}} L(g; \mathbf{m}) \right].$$

*In the homogeneous case, $\rho(du|s) = \rho(du)$, the expression reduces to*

$$\sum_{\mathbf{p}} \left[ \sum_{\mathbf{m} \in S_{n(\mathbf{p})}} \left[ \prod_{j=1}^{n(\mathbf{p})} \kappa_{m_j, r_{j-1}}(\rho) \right] \right.$$

$$\left. \times \int_0^\infty \int_{t_{n(\mathbf{p})}}^\infty \cdots \int_{t_2}^\infty g(t, \mathbf{m}) \prod_{j=1}^{n(\mathbf{p})} e^{-\Lambda_0(t_j)\psi_{m_j, r_{j-1}}} \Lambda_0(dt_j) \right].$$

PROOF.  The result follows from an application of Fubini's theorem and Lemma 5.2, which yields

$$\int_{\mathcal{M}} \mathscr{I}(g) \mathbb{P}(dN|\nu) = \int g(\mathbf{t}) \mathscr{M}(d\mathbf{t}, d\mathbf{x}). \qquad \square$$

REMARK 8.  The case where $g$ may depend also on $\mathbf{X}$ is obvious. It is important to note that Lemma 5.3 may viewed as a generalization of Lo ([32], Lemma 2).

We now use Lemma 5.3 to obtain a simpler description of the distribution of $(\mathbf{T}, \mathbf{X})$, which also yields easily the EPPF formulae and a corresponding distribution on $(\mathbf{m}, \mathbf{p})$. Note again that we do this without resorting to the types of combinatorial arguments used, for instance, in [1] and [32].

PROPOSITION 5.1.  *Let $(\mathbf{T}, \mathbf{X})$ denote the random variables with the exchangeable distribution $\mathscr{M}$ described in Lemma 5.2. Then this distribution may be expressed in terms of a conditional distribution of $\mathbf{T}, \mathbf{X}|\mathbf{m}, \mathbf{p}$ and a distribution of $(\mathbf{m}, \mathbf{p})$ as follows:*

(i) *There exists a marginal distribution of $\mathbf{T}, \mathbf{X}|\mathbf{m}, \mathbf{p}$, given by $\pi(d\mathbf{T}, d\mathbf{X}|\mathbf{m}, \mathbf{p})$ proportional to*

$$\left[ \prod_{j=1}^{n(\mathbf{p})} e^{-\int_0^{T_{(j:n)}} \psi_{m_j, r_{j-1}}(s) \Lambda_0(ds)} \kappa_{m_j, r_{j-1}}(\rho|T_{(j:n)}) \right] \prod_{l=1}^{n(\mathbf{p})} \Lambda_0(dT_{(j:n)}, dX_j^*),$$

*where $T_{(1:n)} > T_{(2:n)} > \cdots > T_{(n(\mathbf{p}):n)}$ denotes the order statistics of the unique values $\mathbf{T}^*$. In the homogeneous case the result reduces to*

$\pi(d\mathbf{T}, d\mathbf{X}|\mathbf{m}, \mathbf{p})$

(10)
$$= \left[\prod_{j=1}^{n(\mathbf{p})} \phi(r_j)\right]\left[\prod_{j=1}^{n(\mathbf{p})} e^{-\psi_{m_j,r_{j-1}}\Lambda_0(T_{(j\,:\,n)})}\right] \prod_{l=1}^{n(\mathbf{p})} \Lambda_0(dT_{(j\,:\,n)}, dX_j^*).$$

In both cases $\prod_{j=1}^{n(\mathbf{p})} P_0(dX_j^*|T_{(j:n)})$ is the conditional distribution of $\mathbf{X}|\mathbf{T},\mathbf{m}$.

(ii) *The distribution of* $(\mathbf{m},\mathbf{p})$, *is described as follows. The EPPF derived by i.i.d. sampling from $F$ is expressible as*

$$\pi(\mathbf{p}) = \sum_{\mathbf{m}\in S_{n(\mathbf{p})}} L(1;\mathbf{m}).$$

*The representations imply the existence of a joint distribution of $(\mathbf{m},\mathbf{p})$ given by $\pi(\mathbf{m},\mathbf{p}) = L(1;\mathbf{m})$. Additionally, in the case where $\rho(du|s) = \rho(du)$, the formulae reduce to*

$$(11)\quad \pi(\mathbf{p}) = \sum_{\mathbf{m}\in S_{n(\mathbf{p})}} \frac{\prod_{j=1}^{n(\mathbf{p})} \kappa_{m_j,r_{j-1}}(\rho)}{\prod_{j=1}^{n(\mathbf{p})} \phi(r_j)} \quad and \quad \pi(\mathbf{m},\mathbf{p}) = \frac{\prod_{j=1}^{n(\mathbf{p})} \kappa_{m_j,r_{j-1}}(\rho)}{\prod_{j=1}^{n(\mathbf{p})} \phi(r_j)}.$$

PROOF. Statement (i) follows from (ii) and Lemma 5.2. The proof of (ii) in the general case follows from Lemma 5.3 with $g := 1$. In the case of $\rho(du|s) = \rho(du)$, $\pi(\mathbf{p})$ is equivalent to

$$\sum_{\mathbf{m}\in S_{n(\mathbf{p})}} \left[\prod_{j=1}^{n(\mathbf{p})} \kappa_{m_j,r_{j-1}}(\rho)\right] \int_0^\infty \int_{t_{n(\mathbf{p})}}^\infty \cdots \int_{t_2}^\infty \prod_{j=1}^{n(\mathbf{p})} e^{-\Lambda_0(t_j)\psi_{m_j,r_{j-1}}} \Lambda_0(dt_j).$$

The result is concluded by evaluating $\int_0^\infty \int_{t_{n(\mathbf{p})}}^\infty \cdots \int_{t_2}^\infty \prod_{j=1}^{n(\mathbf{p})} e^{-\Lambda_0(t_j)\psi_{m_j,r_{j-1}}} \times \Lambda_0(dt_j)$. This is done by noting that for any positive $C$, $\int_t^\infty e^{-C\Lambda_0(u)}\Lambda_0(du) = C^{-1}e^{-C\Lambda_0(t)}$. In addition, $\psi_{m_1,0} = \phi(m_1)$, and for each $j$, $\phi(r_{j-1}) + \psi_{m_j,r_{j-1}} = \phi(r_j)$. □

Equation (10) in Proposition 5.1 can be used to deduce an explicit Markov property in the homogeneous case that has the interpretation that the distribution of the next death time only depends on the previous death time. Moreover, it demonstrates that it is fairly simple to sample from (10).

PROPOSITION 5.2. *Given* $(\mathbf{m},\mathbf{p})$, *let* $T_{(1\,:\,n)},\ldots,T_{(n(\mathbf{p})\,:\,n)}$ *be distributed according to* (10). *Moreover, set* $\Lambda_0(t) = t$. *Then, conditional on* $T_{(j+1\,:\,n)},\ldots,T_{n(\mathbf{p})}$, *the distribution of* $T_{(j\,:\,n)}$ *depends only on* $T_{(j+1)} = t_{j+1}$ *and is given by the truncated exponential distribution with density*

$$\mathbb{P}(T_{(j\,:\,n)} \in dt_j | T_{(j+1\,:\,n)} = t_{j+1}) = \phi(r_j) e^{-\phi(r_j)[t_j - t_{j+1}]} dt_j$$

for $t_j > t_{j+1}$. In particular, the smallest value, or equivalently the first of $n(\mathbf{p})$ death times, $T_{(n(\mathbf{p}):n)}$, has a marginal distribution that is exponential with parameter $\phi(n)$, that is,

$$\mathbb{P}(T_{(n(\mathbf{p}):n)} \in dy) = \phi(n)e^{-\phi(n)y}\,dy.$$

5.1. *Some connections to exponential functionals and means of NTR processes.* We now relate some of our results to those of Epifani, Lijoi and Prünster [11] and Carmona, Petit and Yor [5] concerning moment formulae for means of NTR processes. Briefly, using the relationship in (2), Epifani, Lijoi and Prünster ([11], Proposition 5) established the following moment formulae, expressed in our notation, that characterizes the distribution of $I$:

$$\mathbb{E}[I^n|\nu] = n! \int_0^\infty \int_{t_{n(\mathbf{p})}}^\infty \cdots \int_{t_2}^\infty \prod_{j=1}^n \exp\left(-\int_0^{t_j}\int_0^\infty (1-e^{-y})\right.$$

(12)
$$\left. \times e^{-y(j-1)}\tau(dy|s)\Lambda_0(ds)\right)dt_j.$$

The authors also provide conditions under which the moments exist, which amounts to the finiteness of the moment of order $n$ of $F_0$; that is, $\int_0^\infty t^n F_0(dt) < \infty$. In addition, when $\rho(du|s) = \rho(du)$ and $\Lambda_0(t) = t$, the expression in (12) reduces to the interesting formulae of Carmona, Petit and Yor ([5], Proposition 3.3), viewed within the context of exponential functionals of a subordinator,

$$\mathbb{E}[I^n|\nu] = \frac{n!}{\prod_{j=1}^n \phi(j)}. \tag{13}$$

Notice that the specification $\Lambda_0(t) = t$ is equivalent to specifying $F_0$ as an exponential(1) distribution. In addition, Carmona, Petit and Yor ([5], Proposition 3.1) establish the following result for any $\lambda \geq 1$ and more general Lévy processes:

$$\mathbb{E}[I^\lambda] = \frac{\lambda}{\phi(\lambda)} E[I^{\lambda-1}].$$

Lemma 5.3 offers a complementary result to theirs in that one can express $\mathbb{E}[I^n|\nu]$ in terms of sums over partitions $\mathbf{p}$. Apparently, for NTR processes, a result of this type is only widely known in the case of the Dirichlet process, which follows as a special case of Lo [32]. The result is as follows.

COROLLARY 5.1. *Let $I$ be defined as in* (2). *Then setting* $g(\mathbf{t}) = \prod_{i=1}^n t_i = \prod_{j=1}^{n(\mathbf{p})} t_{(j)}^{m_j}$ *in Lemma* 5.3, *one has* $I = \mathscr{I}(g)$ *and hence*

$$\mathbb{E}[I^n|\nu] = \sum_{\mathbf{p}}\left[\sum_{\mathbf{m}\in S_{n(\mathbf{p})}} L(g;\mathbf{m})\right].$$

*In particular, in the case where $\rho(du|s) = \rho(du)$ and $\Lambda_0(t) = t$, Lemma* 5.3 *combined with the result of Carmona, Petit and Yor* [5] *yields the identity*

$$\sum_{\mathbf{p}} \left[ \sum_{\mathbf{m} \in S_{n(\mathbf{p})}} \left[ \prod_{j=1}^{n(\mathbf{p})} \kappa_{m_j, r_{j-1}}(\rho) \right] \right.$$
$$\left. \times \int_0^\infty \int_{t_{n(\mathbf{p})}}^\infty \cdots \int_{t_2}^\infty \prod_{j=1}^{n(\mathbf{p})} t_j^{m_j} e^{-t_j \psi_{m_j, r_{j-1}}} \, dt_j \right] = \frac{n!}{\prod_{j=1}^n \phi(j)}.$$

Another relationship to the formula for $\mathbb{E}[I^n|\nu]$, (13), given in [5], is seen in the next corollary, derived from Proposition 5.1, which describes the formula for the case where all cells are of the same size.

COROLLARY 5.2. *Suppose that $\rho(du|s) = \rho(du)$ and $n = kn(\mathbf{p})$. Then with respect to the EPPF given in* (11), *the probability of the event $\mathbf{p} = \{E_1, \ldots, E_{n(\mathbf{p})}\}$, such that the size of each cell is $k$, is*

$$\pi(\mathbf{p}) = \frac{n(\mathbf{p})! \prod_{j=1}^{n(\mathbf{p})} \int_0^1 u^k (1-u)^{(j-1)k} \rho(du)}{\prod_{j=1}^{n(\mathbf{p})} \phi(jk)}.$$

*As special cases, when $n(\mathbf{p}) = n$, the probability of no ties in the sample corresponds to the probability of the event $p = \{\{1\}, \{2\}, \ldots, \{n\}\}$ given by*

$$\pi(\mathbf{p}) = \frac{n! \prod_{j=1}^n \int_0^1 u(1-u)^{j-1} \rho(du)}{\prod_{j=1}^n \phi(j)} = \mathbb{E}[I^n|\nu] \prod_{j=1}^n \int_0^1 u(1-u)^{j-1} \rho(du),$$

*for $\mathbb{E}[I^n|\nu]$ given in* (13). *When $n(\mathbf{p}) = 1$, $\mathbf{p} = \{1, 2, \ldots, n\}$ corresponds to the event that all the values in the sample are the same, and the probability is given by*

$$\pi(\mathbf{p}) = \frac{\int_0^1 u^n \rho(du)}{\phi(n)} = \frac{\int_0^\infty (1 - e^{-y})^n \tau(dy)}{\int_0^\infty (1 - e^{-ny}) \tau(dy)}.$$

REMARK 9. The event of no ties, $n(\mathbf{p}) = n$, corresponds to the common assumption in the literature for observed data. Analogous to Antoniak [1] for the Dirichlet process, it follows that when $n(\mathbf{p}) = n$, using Corollary 5.2, the distribution of $\mathbf{T}, \mathbf{X}|\mathbf{p}$ in the homogeneous case is

$$\mathbb{E}[I^n|\nu]^{-1} \left[ \prod_{i=1}^n e^{-\Lambda_0(T_{(i:n)}) \psi_{1,i-1}} \right] \prod_{j=1}^n \Lambda_0(dT_j, dX_j).$$

REMARK 10. Gnedin and Pitman [18], independent of this work and by different arguments, obtain formulae for what are called regenerative compositions that contain our results in (11). Their formulae are derived from a discretization of subordinators. In fact, the authors show that all such regenerative compositions are determined uniquely by their construction via subordinators. The authors' result is more general, in the homogeneous case, because they include the result for subordinators with drift components. It is, however, a simple matter to adjust our results to allow for a drift (see [24], Remark 28). They do not cover the inhomogeneous cases we consider. We discovered these connections through a mutual exchange of manuscripts in progress. The authors' description via a *decrement function* and *composition structure* contain additional binomial coefficients. Explicitly in terms of our notation, their composition structure is expressed as

$$\frac{n!}{\prod_{j=1}^{n(\mathbf{p})} e_j!} \pi(\mathbf{m}, \mathbf{p}).$$

The authors identify some particularly interesting composition structures and we will show how this translates into an interesting class of spatial NTR models. See also [10, 17, 35] for relevant references. See also [19, 20] for important results related to the rates of various $n(\mathbf{p})$.

5.2. *Sampling $\mathcal{M}$: modified Chinese restaurant processes.* Propositions 5.1, 5.2 and 4.1 dictate how one might sample $(\mathbf{T}, \mathbf{X})$ from $\mathcal{M}$. This is especially true in the homogeneous case. One proceeds essentially by first obtaining a draw of $(\mathbf{m}, \mathbf{p})$ from $\pi(\mathbf{m}, \mathbf{p})$, then using Proposition 5.1 or 5.2 to draw the ordered unique values $(T_{(j:n)})$ from the relevant truncated exponential distributions. The $X_j^*$ are then drawn from $P_0(dX_j^*|T_{(j:n)})$ for $j = 1, \ldots, n(\mathbf{p})$. Additionally one can then (approximately) draw $F|\mathbf{T}, \mathbf{X}$, by using the representation $F_n^*$ from Proposition 4.1, which suggests to draw $(J_{j,n})$, and then applying methods in the literature to approximate quantities such as $\Lambda_n$ (see, e.g., [4]). These are precisely the type of steps that would lead to efficient approximations in more complex mixture models, that is to say, models where $(\mathbf{T}, \mathbf{X})$ are missing values obtained from $\mathcal{M}$ and are not directly observed. Also, by sampling from $\mathcal{M}$ one can approximate quantities such as those that appear in Lemma 5.3. In this section it is shown how one might generate $(\mathbf{m}, \mathbf{p})$ from $\pi(\mathbf{m}, \mathbf{p})$ in the case where $\rho(du|s) = \rho(du)$ via a sequential seating scheme with probabilities derived from the prediction rule given $(\mathbf{m}, \mathbf{p})$. This idea also holds in the nonhomogeneous case. The scheme bears similarities to generalized Chinese restaurant processes that can be used to generate general EPPFs, $\pi(\mathbf{p}) = p(e_1, \ldots, e_{n(\mathbf{p})})$. Using the description in [37], page 60, the generalized Chinese restaurant scheme assumes that an initially empty Chinese restaurant has an unlimited number of tables labeled $1, 2, \ldots$. Customers numbered $1, 2, \ldots$ arrive one by and are seated sequentially according to probabilities derived from ratios of the EPPF. Basically customers are seated with probabilities that depend on the size or number of customers already seated at the existing tables.

5.2.1. *Ordered generalized Chinese restaurant processes.* In general, to draw from $\tilde{p}(m_1, \ldots, m_{n(\mathbf{p})}) := \pi(\mathbf{m}, \mathbf{p})$, we introduce a new scheme, which is a modified Chinese restaurant process that also records the rank of the entering customers relative to the already seated customers. The first customer is seated and assigned an initial rank of 1. Now, given a configuration based on $n$ customers seated at $n(\mathbf{p})$ existing tables labeled with ranks from $j = 1, \ldots, n(\mathbf{p})$, the next customer $n+1$ is seated at an occupied table $j$, denoting that customer $n+1$ is equivalent to the $j$th largest seated customers, with probability

$$\begin{aligned}
p_{j:n} &= \frac{\tilde{p}(\ldots, m_j + 1, \ldots)}{\tilde{p}(m_1, \ldots, m_{n(\mathbf{p})})} \\
&= \frac{\kappa_{m_j+1, r_{j-1}}(\rho) \prod_{l=j+1}^{n(\mathbf{p})} \kappa_{m_l, r_{l-1}+1}(\rho)}{\kappa_{m_j, r_{j-1}}(\rho) \prod_{l=j+1}^{n(\mathbf{p})} \kappa_{m_l, r_{l-1}}(\rho)} \prod_{l=j}^{n(\mathbf{p})} \frac{\phi(r_l)}{\phi(r_l + 1)}.
\end{aligned} \tag{14}$$

Customer $n+1$ is seated at a new table with probability $1 - \sum_{j=1}^{n(\mathbf{p})} p_{j:n}$. However, if customer $n+1$ is new, it is also necessary to know the customer's rank and as such to rerank by one position all customers smaller than the new customer. Hence the probability that customer $n+1$ is new and is the $j$th largest among $n(\mathbf{p}) + 1$ possible ranks is

$$\begin{aligned}
q_{j:n} &= \frac{\tilde{p}(\ldots, m_{j-1}, 1, m_j, \ldots)}{\tilde{p}(m_1, \ldots, m_{n(\mathbf{p})})} \\
&= \frac{\kappa_{1, r_{j-1}}(\rho)}{\phi(r_{j-1}+1)} \frac{\prod_{l=j}^{n(\mathbf{p})} \kappa_{m_l, r_{l-1}+1}(\rho)}{\prod_{l=j}^{n(\mathbf{p})} \kappa_{m_l, r_{l-1}}(\rho)} \prod_{l=j}^{n(\mathbf{p})} \frac{\phi(r_l)}{\phi(r_l+1)}
\end{aligned}$$

with $q_{n(\mathbf{p})+1:n} = \kappa_{1,n}(\rho)/\phi(n+1)$. Note that in the calculation of $\kappa_{1, r_{j-1}}(\rho)$, $r_{j-1} + 1$ is to be used rather than $r_j = r_{j-1} + m_j$.

As an example, consider the choice of a homogeneous beta process that corresponds to $c(s) = \theta$ in (4). Then it is easily seen that $\phi(r_j) = \sum_{l=1}^{r_j} \theta/(\theta + l - 1)$ and it follows that, in this case,

$$p_{j:n} = \frac{m_j}{n+\theta} \prod_{l=j}^{n(\mathbf{p})} \frac{\phi(r_l)}{\phi(r_l+1)}$$

and

$$q_{j:n} = \frac{1}{n+\theta} \frac{1}{\sum_{i=1}^{r_{j-1}+1} 1/(\theta+i-1)} \prod_{l=j}^{n(\mathbf{p})} \frac{\phi(r_l)}{\phi(r_l+1)}.$$

5.3. *Species sampling models generated by spatial NTR processes.* The availability of the EPPF, coupled with Pitman's [34] theory of species sampling random probability models, implies that there exists a new explicit class of random probability measures of the form

$$P_F(\cdot) = \int_0^\infty S(s-)\Lambda(ds,\cdot) = \sum_{i=1}^\infty Q_i \delta_{Z_i}(\cdot), \tag{15}$$

where $Z_i$ are i.i.d. random elements in $\mathscr{X}$ with some nonatomic law $P_0$ and where, independent of $(Z_i)$, $(Q_i)$ denotes a collection of random probabilities that sum to 1 and whose law is completely determined by the EPPF $\pi(\mathbf{p})$ given in Proposition 5.1. We will call $P_F$ an *NTR species sampling model*. We do point out that although there are technically a large number of possible species sampling models, to date there are only two well-known classes: the species sampling models based on the Poisson–Kingman models described in [38] (see also [24]) and the stick-breaking models described in [22]. See also [40].

The NTR species sampling model, which is defined for the first time here, represents a third case where, due to the present analysis, much is known. All three classes contain the Dirichlet process. In fact, rather remarkably, all three classes contain the two-parameter $(\alpha, \theta)$ Poisson–Dirichlet family of random probability measures with parameters $0 \leq \alpha < 1$ and $\theta \geq 0$. We will describe this in a forthcoming section. The next proposition describes how one can always formally obtain an NTR species sampling model generated by an $F$ with an independent prior specification, $F_0(ds, dx) = F_0(ds)P_0(dx)$. Moreover, we give a description of its posterior distribution.

PROPOSITION 5.3. *Let $\nu(du, ds, dx) = \rho(du|s)\Lambda_0(ds, dx)$ denote the mean intensity of a Poisson random measure $N$ on $\mathscr{W}$, where $\Lambda_0$ is chosen such that $\Lambda_0(ds, dx) = \Lambda_0(ds)P_0(dx)$. Then the corresponding spatial NTR process, $F$, generates an NTR species sampling model, $P_F$, given in* (15), *by the representations $P_F(dx) := F(\infty, dx) = \int_0^\infty S(s-)\Lambda(ds, dx)$ or, equivalently, the marginal distribution of $\mathbf{X} = (\mathbf{X}^*, \mathbf{p})$ is given by*

$$\mathbb{E}\left[\prod_{i=1}^n P_F(dX_i)|\nu\right] = \pi(\mathbf{p}) \prod_{j=1}^{n(\mathbf{p})} P_0(dX_j^*).$$

*Additionally, the posterior distribution of $P_F$ given $(\mathbf{T}, \mathbf{X})$, or just $\mathbf{X}$, is characterized by Proposition* 4.1 *and* 5.1. *Specifically, it is equivalent to the appropriate conditional laws of the random measure $F_n^*(\infty, dx)$.*

PROOF. Under the specifications $F_0(ds, dx) = F_0(ds)P_0(dx)$, $\mathscr{M}(d\mathbf{T}, d\mathbf{X})$ is such that given $\mathbf{p}$, the vectors $\mathbf{T}^*$ and $\mathbf{X}^*$ are independent, where $\mathbf{X}^*$ has joint law $\prod_{j=1}^{n(\mathbf{p})} P_0(dX_j^*)$. The result is concluded by integrating out $\mathbf{T}^*$. □

It is interesting to note that while the Dirichlet process is an example of $P_F$, it also arises without the independence specification. In most cases $\pi(\mathbf{p})$ will not be easy to work with directly; as such one can work with $\pi(\mathbf{m}, \mathbf{p})$. As an example, we present a description for the prediction rule of $P_F$ given $(\mathbf{X}, \mathbf{m})$. It will be clear that one can employ the *ordered generalized Chinese restaurant* algorithm in Section 5.2.1 to draw easily from a joint distribution of $(\mathbf{X}, \mathbf{m})$.

PROPOSITION 5.4. *Let $P_F$ denote an NTR species sampling model defined by the choice $\rho(du|s) = \rho(du)$. Suppose that $\mathbf{X} = \{X_1, \ldots, X_n\}$ given $P_F$ are i.i.d. $P_F$. Then one can define a prediction rule for $X_{n+1}$ given $\mathbf{X}, \mathbf{m}$ as*

$$\mathbb{P}(X_{n+1} \in dx | \mathbf{X}, \mathbf{m}) = \left(1 - \sum_{j=1}^{n(\mathbf{p})} p_{j\,:\,n}\right) P_0(dx) + \sum_{j=1}^{n(\mathbf{p})} p_{j\,:\,n} \delta_{X_j^*}(dx),$$

*where $(p_{j\,:\,n})$ are given in (14). Note also that $\pi(\mathbf{m}, \mathbf{p}) \prod_{i=1}^{n(\mathbf{p})} P_0(dX_j)$ is the distribution of $(\mathbf{X}, \mathbf{m})$, which means that the distribution of $\mathbf{X}|\mathbf{m}$ is such that the unique values $(X_j^*)$ given $\mathbf{m}$ are i.i.d. $P_0$. The prediction rule given $\mathbf{X}$ is obtained by $\mathbb{P}(X_{n+1} \in dx|\mathbf{X}) = \sum_{\mathbf{m} \in S_{n(\mathbf{p})}} \mathbb{P}(X_{n+1} \in dx|\mathbf{X}, \mathbf{m})\pi(\mathbf{m}|\mathbf{p})$.*

## 6. Examples.

6.1. *Generalized gamma models.* An interesting class of measures is the family of generalized gamma random measures discussed in [4]. Using the description of Brix [4], these are $Z$ processes with Lévy measure

$$\tau_{\alpha,b}(dy)\Lambda_0(ds, dx) = \frac{1}{\phi_{\alpha,b}(1)\Gamma(1-\alpha)} y^{-\alpha-1} \exp(-by)\, dy \Lambda_0(ds, dx),$$

where $\phi_{\alpha,b}(1) = \frac{1}{\alpha}[(b+1)^\alpha - (b)^\alpha]$. The values for $\alpha$ and $b$ are restricted to satisfy $0 < \alpha < 1$ and $0 \leq b < \infty$ or $-\infty < \alpha \leq 0$ and $0 < b < \infty$. Different choices for $\alpha$ and $b$ in $\rho_{\alpha,b}$ yield various subordinators. These include the stable subordinator when $b = 0$, the gamma process subordinator when $\alpha = 0$ and the inverse-Gaussian subordinator when $\alpha = 1/2$ and $b > 0$. When $\alpha < 0$, this results in a class of gamma compound Poisson processes. Generalized gamma NTR processes with $b > 0$ are discussed in [11]. Here, from our results,

$$\psi_{m_j, r_{j-1}} = \frac{[(r_j + b)^\alpha - (r_{j-1} + b)^\alpha]}{[(1+b)^\alpha - b^\alpha]},$$

$$\phi(r_j) = \frac{[(r_j + b)^\alpha - b^\alpha]}{[(1+b)^\alpha - b^\alpha]}$$

and

$$\kappa_{m_j,r_{j-1}}(\rho) = \frac{\sum_{l=0}^{m_j}(-1)^{l+1}\binom{m_j}{l}(b+r_{j-1}+l)^\alpha}{[(1+b)^\alpha - b^\alpha]}.$$

Hence

$$\pi(\mathbf{m},\mathbf{p}) = \frac{\prod_{j=1}^{n(\mathbf{p})}\sum_{l=0}^{m_j}\binom{m_j}{l}(-1)^{l+1}(b+r_{j-1}+l)^\alpha}{\prod_{j=1}^{n(\mathbf{p})}[(b+r_j)^\alpha - b^\alpha]}.$$

The process $F(ds,dx)$ is such that marginally $F(ds,\mathscr{X})$ is a generalized gamma NTR process and $P_F(dx) = F(\infty,dx)$ is a species sampling model. Additionally one can use Proposition 5.2 to generate the $T_{(j:n)}$. In particular, when $b = 0$ and $\Lambda_0(t) = t$, the density corresponding to the stable process with index $0 < \alpha < 1$ is

$$\mathbb{P}(T_{(j:n)} \in dt_j | T_{(j+1:n)} = t_{j+1}) = r_j^\alpha e^{-r_j^\alpha[t_j - t_{j+1}]} \qquad \text{for } t_j > t_{j+1}.$$

6.2. *The spatial NTR two-parameter Poisson–Dirichlet model.* We now describe perhaps the most remarkable class of spatial NTR processes. Gnedin and Pitman ([18], Section 10) were able to deduce that one can generate the EPPF of the two-parameter $(\alpha,\theta)$ Poisson–Dirichlet distribution with parameters $0 \leq \alpha < 1$ and $\theta \geq 0$ by specifying a homogeneous $\rho$ such that

(16) $$\int_u^1 \rho(dv) = \frac{\Gamma(\theta + 2 - \alpha)}{\Gamma(1-\alpha)\Gamma(1+\theta)} u^{-\alpha}(1-u)^\theta$$

and, hence,

$$\phi(r_j) = \frac{r_j \Gamma(\theta + r_j)\Gamma(\theta + 2 - \alpha)}{\Gamma(\theta + 1)\Gamma(\theta - \alpha + r_j + 1)}.$$

Due to Proposition 5.2, this is enough to generate the distribution of the $(T_{(j:n)})$. Note for this model one can directly sample from the well-known EPPF.

6.2.1. *The ordered ESF and the Dirichlet process.* An interesting case is when $\alpha = 0$; that is, $\rho(du) = \theta(\theta+1)(1-u)^{\theta-1} du$. This choice generates the ordered Ewens sampling formula as described in [10]. Moreover, the spatial NTR process $F(ds,dx)$ is such that $F(ds,\mathscr{X})$ is an NTR process but not a Dirichlet process, and it follows from Proposition 5.3 that $P_F(dx) = F(\infty,dx)$ is a Dirichlet process with shape $\theta P_0$. Hence, this shows that a Dirichlet process may be generated via a homogeneous NTR process derived from a compound Poisson process. Note of course that when $x \in \mathcal{R}$, this process is marginally an NTR process in both coordinates. Setting $\Lambda_0(t) = t$, the corresponding distribution of the $T_{(j:n)}$ is given by

$$\mathbb{P}(T_{(j:n)} \in dt_j | T_{(j+1:n)} = t_{j+1}) = \frac{(\theta+1)r_j}{(\theta+r_j)} e^{-((\theta+1)r_j/(\theta+r_j))[t_j - t_{j+1}]}$$

for $t_j > t_{j+1}$.

Note also that the distribution of the jumps that depends on $(T_{(j:n)})$ is

$$\mathbb{P}(J_{j,n} \in du | T_{(j:n)}) = \frac{\Gamma(\theta + r_j + 1)}{\Gamma(m_j + 1)\Gamma(\theta + r_{j-1})} u^{m_j}(1-u)^{\theta + r_{j-1} - 1};$$

that is, they are beta distributed with parameters $(m_j + 1, \theta + r_{j-1})$. Note that these are not the jumps of a posterior Dirichlet process. However, since $F(\infty, dx)$ is a Dirichlet process, its posterior distribution given $\mathbf{X}$ is a Dirichlet process with shape $\theta P_0 + \sum_{i=1}^n \delta_{X_i}$.

6.2.2. *Representations for the general two-parameter $(\alpha, \theta)$ case.* We can use the result above to provide some new results related to the two-parameter $(\alpha, \theta)$ Poisson–Dirichlet family and spatial NTR processes. For clarity, we first recall the definition of the two-parameter $(\alpha, \theta)$ Poisson–Dirichlet class of random probability measures. The two-parameter $(\alpha, \theta)$ Poisson–Dirichlet random probability measure with parameters $0 \leq \alpha < 1$ and $\theta \geq 0$ has the known representation

$$P_{\alpha,\theta}(dx) = \frac{\mu_{\alpha,\theta}(dx)}{T_{\alpha,\theta}},$$

where $\mu_{\alpha,\theta}$ is a finite random measure on $\mathscr{X}$ with law $\mathbb{P}(d\mu_{\alpha,\theta})$ and where $T_{\alpha,\theta} = \mu_{\alpha,\theta}(\mathscr{X})$ is a random variable. The law of the random measure $\mu_{\alpha,\theta}$ can be described as follows. When $\alpha = 0$, $\mu_{0,\theta}$ is a gamma process with shape $\theta P_0$; hence, $P_{0,\theta}$ is a Dirichlet process with shape $\theta P_0$. When $\theta = 0$, $\mu_{\alpha,0}$ is a stable random measure of index $0 < \alpha < 1$. Note that both $\mu_{\alpha,0}$ and $\mu_{0,\theta}$ are completely random measures and can be represented in terms of a Poisson random measure. However, this is not true for the case where both $\alpha$ and $\theta$ are positive. Here for $0 < \alpha < 1$ and $\theta > 0$ one has the absolute continuity relationship

$$\mathbb{P}(d\mu_{\alpha,\theta}) = \frac{T_{\alpha,0}^{-\theta} \mathbb{P}(d\mu_{\alpha,0})}{\mathbb{E}[T_{\alpha,0}^{-\theta}]},$$

where $T_{\alpha,0}$ is a stable law random variable. This class of models also has a representation in terms of stick-breaking processes. See, for instance, [22, 34, 39] for further details. We now arrive at the following interesting observations.

PROPOSITION 6.1. *Let $F(ds, dx)$ denote a spatial NTR process specified by the choice of $\rho$ in* (16) *and let $F_0(ds, dx) = P_0(dx)F_0(ds)$. Then $P_F$ is a two-parameter $(\alpha, \theta)$ Poisson–Dirichlet process. This yields the representations*

$$P_F(dx) = \int_0^\infty S(s-)\Lambda(ds, dx)$$
$$= \sum_{k=1}^\infty V_k \prod_{i=1}^{k-1}(1 - V_i)\delta_{Z_k}(dx)$$
$$= \mu_{\alpha,\theta}(dx)/T_{\alpha,\theta} = P_{\alpha,\theta}(dx),$$

*where $(V_k)$ are independent beta $(1-\alpha, \theta+k\alpha)$ random variables independent of the $(Z_k)$, which are i.i.d. $P_0$; that is, a two-parameter $(\alpha, \theta)$ Poisson–Dirichlet process can be represented as the marginal probability measure of a spatial NTR process as described above.*

PROOF. The general result follows from an application of Proposition 5.3 combined with the calculations of the EPPF using $\rho$ in (16) by Gnedin and Pitman ([18], Section 10). The case of the Dirichlet process that corresponds to the choice of $\rho(du) = \theta(\theta+1)(1-u)^{\theta-1} du$ could be deduced as well from [10] in combination with Proposition 5.3. See also [35] for the $(\alpha, \alpha)$ model. □

## APPENDIX

**Proofs of Proposition 4.1 and Lemmas 5.1 and 5.2.** We now show that the proofs of Proposition 4.1 and Lemmas 5.1 and 5.2 follow as a simple consequence of the Poisson partition calculus methods as laid out in [24, 26]. First set $W_i = (J_i, T_i, X_i)$ for $i = 1, \ldots, n$, elements of $\mathscr{W}$. The collection $\mathbf{J} = \{J_1, \ldots, J_n\}$ with values in $[0, 1]$ will play the role of the latent jumps. Its unique values are the $(J_{j,n})$. Set $\mathbf{W} = (\mathbf{J}, \mathbf{T}, \mathbf{X})$ and let $W_j^* = (J_{j,n}, T_{(j:n)}, X_j^*)$ denote the $j = 1, \ldots, n(\mathbf{p})$ unique triples. Using Proposition 2.3 of [26] yields the following statement. Suppose that $(\mathbf{W}, N)$ are measurable elements in the space $\mathscr{W}^n \times \mathbb{M}$, where $N$ is Poisson random measure with sigma finite nonatomic mean measure $\nu$. Then for each nonnegative measurable $f$ such that $\mathscr{G}(f) < \infty$, the following disintegration holds:

(17)
$$\left[\prod_{i=1}^n N(dW_i)\right] e^{-N(f)} \mathbb{P}(dN|\nu)$$
$$= e^{-\mathscr{G}(f)} \mathbb{P}(dN|\nu_f, \mathbf{W}) \prod_{j=1}^{n(\mathbf{p})} e^{-f(W_j^*)} \nu(dW_j^*),$$

where $\mathbb{P}(dN|\nu_f, \mathbf{W})$ denotes the law of the random measure $N + \sum_{j=1}^{n(\mathbf{p})} \delta_{W_j^*}$, where $N$ is a Poisson random measure with mean intensity $\mathbb{E}[N(du, ds, dx)|\nu_f] = \nu_f(du, ds, dx) := e^{-f(u,s,x)}\nu(du, ds, dx)$.

To apply the results above we first express (9) in terms of an exponential functional of a Poisson random measure as follows. For each $j$, set $f_{T_{(j:n)-}}(u, s, x) = I\{s < T_{(j:n)}\}[-\log(1 - u)]$. Now it follows that one can define

$$f_n(u, s, x) = \sum_{j=1}^{n(\mathbf{p})} m_j f_{T_{(j:n)-}}(u, s, x) = -Y_n(s)\log(1 - u)$$

and hence one has

$$\left[\prod_{j=1}^{n(\mathbf{p})} S(T_{(j:n)}-)^{m_j}\right] = \prod_{j=1}^{n(\mathbf{p})} e^{-N(m_j f_{T_{(j:n)-}})} = e^{-N(f_n)}.$$

Note also that $e^{-f_n(u,s,x)} = (1-u)^{Y_n(s)}$ and $e^{-f_n(J_{j,n}, T_{(j:n)}, X_j^*)} = (1 - J_{j,n})^{r_j-1}$ for $j = 1, \ldots, n(\mathbf{p})$.

The next step is to write $\Lambda(dT_i, dX_i) = \int_0^1 J_i N(dJ_i, dT_i, dX_i)$. Now removing those integrals in (6) yields an augmentation of the distribution of $(\mathbf{T}, \mathbf{X}, N)$ in terms of a distribution of $(\mathbf{J}, \mathbf{T}, \mathbf{X}, N)$. It follows that the distribution of $(\mathbf{J}, \mathbf{T}, \mathbf{X}, N)$ can be expressed similar to the left-hand side of (17) with $f_n$ in place of $f$ as

$$\left[\prod_{j=1}^{n(\mathbf{p})} J_{j,n}^{m_j}\right]\left[\prod_{i=1}^{n} N(dJ_i, dT_i, dX_i)\right] e^{-N(f_n)}\mathbb{P}(dN|\nu).$$

Note that $\prod_{i=1}^{n} J_i = \prod_{j=1}^{n(\mathbf{p})} J_{j,n}^{m_j}$. Hence now applying the right-hand side of (17) one has that the joint distribution of $(\mathbf{J}, \mathbf{T}, \mathbf{X}, N)$ is given by

(18)
$$\mathbb{P}(dN|e^{-f_n}\nu_{f_n}, \mathbf{W})e^{-\mathscr{G}(f_n)}$$
$$\times \left[\prod_{j=1}^{n(\mathbf{p})} J_{j,n}^{m_j}(1 - J_{j,n})^{r_j-1}\rho(dJ_{j,n}|T_{(j:n)})\right]\prod_{l=1}^{n(\mathbf{p})} \Lambda_0(T_l^*, X_l^*),$$

where $\mathbb{E}[e^{-N(f_n)}|\nu] = e^{-\mathscr{G}(f_n)}$ and now $\mathbb{P}(dN|\nu_{f_n}, \mathbf{W})$ corresponds to the law, for fixed $\mathbf{W}$, of a random measure $N_n + \sum_{j=1}^{n(\mathbf{p})} \delta_{J_{j,n}, T_{(j:n)}, X_j^*}$, where $N_n$ is a Poisson random measure with mean described in (8) and $\mathbb{P}(dN|\nu_{f_n}, \mathbf{W})$ is the posterior distribution of $N|\mathbf{J}, \mathbf{T}, \mathbf{X}$. The joint distribution of $(\mathbf{J}, \mathbf{T}, \mathbf{X})$ is obtained by integrating out $N$ in (18). Now using the fact that one can decompose $(\mathbf{J}, \mathbf{T}, \mathbf{X})$ as $((J_{j,n}), \mathbf{T}^*, \mathbf{X}^*, \mathbf{p})$, it follows that an expression for the marginal distribution of $(\mathbf{T}, \mathbf{X})$, or equivalently $(\mathbf{T}^*, \mathbf{X}^*, \mathbf{p})$, is obtained by integrating out $N$ and the $(J_{j,n})$. For clarity, this takes the form

$$\mathscr{M}(d\mathbf{T}, d\mathbf{X}) = e^{-\mathscr{G}(f_n)} \left[ \prod_{j=1}^{n(\mathbf{p})} \kappa_{m_j, r_{j-1}}(\rho | T_{(j:n)}) \right] \prod_{l=1}^{n(\mathbf{p})} \Lambda_0(dT_l^*, X_l^*).$$

The description of the posterior distribution of $N|\mathbf{T}, \mathbf{X}$ is given in terms of the distribution of $N|\mathbf{J}, \mathbf{T}, \mathbf{X}$ mixed over the distribution of the $(J_{j,n})$ given $(\mathbf{T}, \mathbf{X})$. The distribution of $(J_{j,n})$ follows by an appeal to the classical Bayes rule; that is, one integrates out $N$ in (18) and then divides the remaining quantity by $\mathscr{M}$. This yields the results in Proposition 4.1. Now it follows that the description of $\mathscr{M}$ given in Lemma 5.2 is completed by verifying Lemma 5.1. This is obtained by using repeatedly the exponential change of measure described in Proposition 2.1 of [26]. This is the same as working with (17) after removing all the terms that involve $\mathbf{W}$; that is, the disintegration $e^{-N(f)}\mathbb{P}(dN|\nu) = \mathbb{P}(dN|\nu_f)\mathbb{E}[e^{-N(f)}|\nu]$. We apply this repeatedly to the measure $[\prod_{j=1}^{n(\mathbf{p})} e^{-N(m_j f_{T_{(j:n)-}})}]\mathbb{P}(dN|\nu)$. To see this, first set $g_j := m_j f_{T_{(j:n)-}}$ for $j = 1, \ldots, n(\mathbf{p})$ and let each $g_j$ now play the role of an $f$. We demonstrate the first two steps. Notice that the first term is obtained as

$$e^{-N(g_1)}\mathbb{P}(dN|\nu) = \mathbb{P}(dN|\nu_{g_1})e^{-\int_0^{T_{(1:n)}} \psi_{m_1, r_0}(s) \Lambda_0(ds)}$$
$$= \mathbb{P}(dN|\nu_{g_1})\mathbb{E}[e^{-N(g_1)}|\nu].$$

The next term is obtained as

$$e^{-N(g_2)}\mathbb{P}(dN|\nu_{g_1}) = \mathbb{P}(dN|\nu_{g_1+g_2})e^{-\int_0^{T_{(2:n)}} \psi_{m_2, r_1}(s)\Lambda_0(ds)}.$$

The last expression follows from the fact that for $s < T_{(2:n)}$, $e^{-g_2(s,u,x)} = (1-u)^{m_2}$ and $e^{-g_1(s,u,x)} = (1-u)^{m_1}$ with $r_1 = m_1$. The next term would then exploit this type of relationship for $g_1$, $g_2$ and $g_3$ on $s < T_{(3:n)}$, where $m_1 + m_2 = r_2$. It is clear that continuing in this way leads to the conclusion of Lemma 5.1.

REMARK 11. More details, including an analysis of semiparametric models subject to censoring mechanisms, is given in an older version of this manuscript [25].

**Acknowledgments.** I would like to thank Kjell Doksum for early comments that helped in the exposition of this work. Thanks to Jim Pitman for clarifying some nice connections to his and Alexander Gnedin's work.

# REFERENCES

[1] ANTONIAK, C. E. (1974). Mixtures of Dirichlet processes with applications to Bayesian nonparametric problems. *Ann. Statist.* **2** 1152–1174. MR0365969

[2] BERTOIN, J. and YOR, M. (2001). On subordinators, self-similar Markov processes and some factorizations of the exponential variable. *Electron. Comm. Probab.* **6** 95–106. MR1871698

[3] BLACKWELL, D. and MACQUEEN, J. B. (1973). Ferguson distributions via Pólya urn schemes. *Ann. Statist.* **1** 353–355. MR0362614

[4] BRIX, A. (1999). Generalized gamma measures and shot-noise Cox processes. *Adv. in Appl. Probab.* **31** 929–953. MR1747450

[5] CARMONA, P., PETIT, F. and YOR, M. (1997). On the distribution and asymptotic results for exponential functionals of Lévy processes. In *Exponential Functionals and Principal Values Related to Brownian Motion* (M. Yor, ed.) 73–130. Biblioteca de la Revista Matematica Iberoamericana, Madrid. MR1648657

[6] DALEY, D. J. and VERE-JONES, D. (1988). *An Introduction to the Theory of Point Processes.* Springer, New York. MR0950166

[7] DEY, J. (1999). Some properties and characterizations of neutral-to-the-right priors and beta processes. Ph.D. dissertation, Michigan State Univ.

[8] DEY, J., ERICKSON, R. V. and RAMAMOORTHI, R. V. (2003). Some aspects of neutral to right priors. *Internat. Statist. Rev.* **71** 383–401.

[9] DOKSUM, K. A. (1974). Tailfree and neutral random probabilities and their posterior distributions. *Ann. Probab.* **2** 183–201. MR0373081

[10] DONNELLY, P. and JOYCE, P. (1991). Consistent ordered sampling distributions: Characterization and convergence. *Adv. in Appl. Probab.* **23** 229–258. MR1104078

[11] EPIFANI, I., LIJOI, A. and PRÜNSTER, I. (2003). Exponential functionals and means of neutral-to-the-right priors. *Biometrika* **90** 791–808. MR2024758

[12] EWENS, W. J. (1972). The sampling theory of selectively neutral alleles. *Theoret. Population Biology* **3** 87–112. MR0325177

[13] FERGUSON, T. S. (1973). A Bayesian analysis of some nonparametric problems. *Ann. Statist.* **1** 209–230. MR0350949

[14] FERGUSON, T. S. (1974). Prior distributions on spaces of probability measures. *Ann. Statist.* **2** 615–629. MR0438568

[15] FERGUSON, T. S. and PHADIA, E. (1979). Bayesian nonparametric estimation based on censored data. *Ann. Statist.* **7** 163–186. MR0515691

[16] FREEDMAN, D. A. (1963). On the asymptotic behavior of Bayes estimates in the discrete case. *Ann. Math. Statist.* **34** 1386–1403. MR0158483

[17] GNEDIN, A. V. (1997). The representation of composition structures. *Ann. Probab.* **25** 1437–1450. MR1457625

[18] GNEDIN, A. V. and PITMAN, J. (2005). Regenerative composition structures. *Ann. Probab.* **33** 445–479. MR2122798

[19] GNEDIN, A. V., PITMAN, J. and YOR, M. (2006). Asymptotic laws for regenerative compositions: Gamma subordinators and the like. *Probab. Theory Related Fields.* To appear.

[20] GNEDIN, A. V., PITMAN, J. and YOR, M. (2006). Asymptotic laws for compositions derived from transformed subordinators. *Ann Probab.* **34**. To appear.

[21] HJORT, N. L. (1990). Nonparametric Bayes estimators based on beta processes in models for life history data. *Ann. Statist.* **18** 1259–1294. MR1062708

[22] ISHWARAN, H. and JAMES, L. F. (2001). Gibbs sampling methods for stick-breaking priors. *J. Amer. Statist. Assoc.* **96** 161–173. MR1952729

[23] ISHWARAN, H. and JAMES, L. F. (2003). Generalized weighted Chinese restaurant processes for species sampling mixture models. *Statist. Sinica* **13** 1211–1235. MR2026070

[24] JAMES, L. F. (2002). Poisson process partition calculus with applications to exchangeable models and Bayesian nonparametrics. Available at arxiv.org/abs/math.pr/0205093.

[25] JAMES, L. F. (2003). Poisson calculus for spatial neutral to the right processes. Available at ihome.ust.hk/˜lancelot.

[26] JAMES, L. F. (2005). Bayesian Poisson process partition calculus with an application to Bayesian Lévy moving averages. *Ann. Statist.* **33** 1771–1799. MR2166562

[27] JAMES, L. F. and KWON, S. (2000). A Bayesian nonparametric approach for the joint distribution of survival time and mark variables under univariate censoring. Technical report, Dept. Mathematical Sciences, Johns Hopkins Univ.

[28] KALLENBERG, O. (2002). *Foundations of Modern Probability*, 2nd ed. Springer, New York. MR1876169

[29] KIM, Y. (1999). Nonparametric Bayesian estimators for counting processes. *Ann. Statist.* **27** 562–588. MR1714717

[30] KINGMAN, J. F. C. (1993). *Poisson Processes*. Oxford Univ. Press, New York. MR1207584

[31] LAST, G. and BRANDT, A. (1995). *Marked Point Processes on the Real Line*: *The Dynamic Approach*. Springer, New York. MR1353912

[32] LO, A. Y. (1984). On a class of Bayesian nonparametric estimates. I. Density estimates. *Ann. Statist.* **12** 351–357. MR0733519

[33] LO, A. Y. (1993). A Bayesian bootstrap for censored data. *Ann. Statist.* **21** 100–123. MR1212168

[34] PITMAN, J. (1996). Some developments of the Blackwell–MacQueen urn scheme. In *Statistics*, *Probability and Game Theory* (T. S. Ferguson, L. S. Shapley and J. B. MacQueen, eds.) 245–267. IMS, Hayward, CA. MR1481784

[35] PITMAN, J. (1997). Partition structures derived from Brownian motion and stable subordinators. *Bernoulli* **3** 79–96. MR1466546

[36] PITMAN, J. (1999). Coalescents with multiple collisions. *Ann. Probab.* **27** 1870–1902. MR1742892

[37] PITMAN, J. (2002). Combinatorial stochastic processes. Technical Report 621, Dept. Statistics, Univ. California, Berkeley. Available at stat-www.berkeley.edu/users/pitman/bibliog.html.

[38] PITMAN, J. (2003). Poisson–Kingman partitions. In *Statistics and Science*: *A Festschrift for Terry Speed* (D. R. Goldstein, ed.) 1–34. IMS, Beachwood, OH. MR2004330

[39] PITMAN, J. and YOR, M. (1997). The two-parameter Poisson–Dirichlet distribution derived from a stable subordinator. *Ann. Probab.* **25** 855–900. MR1434129

[40] REGAZZINI, E., LIJOI, A. and PRÜNSTER, I. (2003). Distributional results for means of normalized random measures with independent increments. *Ann. Statist.* **31** 560–585. MR1983542

[41] WALKER, S. and MULIERE, P. (1997). Beta-Stacy processes and a generalization of the Pólya-urn scheme. *Ann. Statist.* **25** 1762–1780. MR1463574

Department of Information
and Systems Management
Hong Kong University
of Science and Technology
Clear Water Bay, Kowloon
Hong Kong
E-mail: lancelot@ust.hk